\documentclass[a4paper,10pt]{article}
\usepackage[latin1]{inputenc}
\usepackage{amsfonts}
\usepackage{amsmath}
\usepackage[all]{xy}
\usepackage{amssymb}
\usepackage{mathabx}
\usepackage{amsthm}
\usepackage[american]{babel}
\usepackage{caption}
\usepackage[OT1,T1]{fontenc}
\usepackage{makeidx}
\usepackage{mathptmx}
\usepackage{multicol}
\usepackage{tikz}
\usepackage{rotating}
\usepackage{subfigure}
\usepackage{upgreek}
\usepackage{xcolor}
\usepackage{lmodern}
\usepackage{graphicx}
\usepackage{pstricks-add}
\usepackage{graphicx}
\usepackage{pst-fun}
\usepackage{psfrag}
\usepackage{ifpdf}
\usepackage{indentfirst}
\usepackage{float}
\usepackage{pgf,tikz,pgfplots}
\pgfplotsset{compat=1.15}
\usepackage{mathrsfs}
\usetikzlibrary{arrows}
\usepackage{hyperref}
\usepackage[title]{appendix}
\usepackage{anyfontsize}

\usepackage{pst-plot}
\usepackage{pst-text}
\newgray{grayteste}{.9}
\usepackage{caption}
\newcommand\raisepunct[1]{\,\mathpunct{\raisebox{0.5ex}{#1}}}

\newcommand{\R}{\mathbb{R}}

\definecolor{wwwwww}{rgb}{0.4,0.4,0.4}
\definecolor{zzttff}{rgb}{0.6,0.2,1}
\definecolor{qqqqff}{rgb}{0,0,1}

\def\m{\mathcal{M}}

\def\dist{\text{dist}}
\def\inte{\text{int }}

\newtheorem{theorem}{\bf Theorem}[section]

\newtheorem{lemma}[theorem]{\bf Lemma}
\newtheorem{corollary}[theorem]{\bf Corollary}

\newtheorem{remark}[theorem]{\bf Remark}

\newtheorem{definition}[theorem]{\bf Definition}

\theoremstyle{definition}
\newtheorem{example}[theorem]{\bf Example}
{\em }

\date{\today}

\begin{document}

\title{Stochastic half-space theorems for minimal surfaces and $H$-surfaces of $\mathbb{R}^{3}$}
\author{G. P. Bessa  \and L. P. Jorge  \and Leandro F. Pessoa }


\date{}

\maketitle


\begin{abstract}

\noindent We prove a version of the strong half-space theorem between the classes of recurrent minimal surfaces and complete minimal surfaces with bounded curvature of $\mathbb{R}^{3}_{\raisepunct{.}}$ We also show that any  minimal hypersurface immersed with bounded curvature in $M\times \R_+$ equals some $M\times \{s\}$ provided $M$ is a complete, recurrent $n$-dimensional Riemannian manifold with $\text{Ric}_M \geq 0$ and whose sectional curvatures are bounded from above. For $H$-surfaces  we prove that a stochastically complete surface $M$ can not be in the mean convex side of a $H$-surface $N$ embedded in $\R^3$ with bounded curvature if $\sup \vert H_{_M}\vert < H$, or ${\rm dist}(M,N)=0$ when $\sup \vert H_{_M}\vert = H$. Finally, a maximum principle at infinity is shown assuming $M$ has non-empty boundary.
\end{abstract}

%

\section{Introduction} 
A classical theorem  in the global theory of minimal surfaces, by  Xavier \cite{xavier}, states that the convex hull of a complete  non-planar  minimal surface of $\mathbb{R}^{3}$ with bounded curvature is the entire $\mathbb{R}^{3}_{\raisepunct{.}}$ This implies, in particular, that  the class of complete minimal surface with bounded curvature has the half-space property, meaning that any complete minimal surface with bounded curvature can not lie in a half-space defined by  some plane $\mathcal{P}\subset \mathbb{R}^{3}$ unless it is a plane parallel to $\mathcal{P}$. In order to show that the  examples of complete minimal surfaces between two parallel planes constructed in \cite{jorge-xavier-ann,rosenberg-toubiana} were not proper,  Hoffman and  Meeks 
in  \cite{hoffmann-meeks},  showed that the class of  properly immersed minimal surfaces of $\mathbb{R}^{3}$  has the half-space property. This  result together with \cite[Thm.8 \& Cor.1]{Meeks-Simon-Yau} yielded the Hoffman-Meeks strong half-space theorem which states that two properly immersed minimal surfaces of $\mathbb{R}^{3}$ intersect unless they are parallel planes.

Likewise, Xavier's half-space theorem yielded a strong half-space theorem for minimal surfaces with bounded curvature, i.e., two complete  minimal surfaces of $\mathbb{R}^{3}$ with bounded curvature must intersect unless they are parallel planes  \cite{bessa-jorge-oliveira,rosenberg-hst}. The proof given in \cite{bessa-jorge-oliveira} yields, as a corollary, a strong half-space theorem between the classes of
  complete proper minimal surfaces and  complete minimal surfaces with bounded curvature of $\mathbb{R}_{\raisepunct{,}}^{3}$ see \cite[Cor.1.4]{bessa-jorge-oliveira}. 
  
It is  worthy mentioning that Neel gave proofs of the Hoffman-Meeks and Xavier's half-space theorems \cite[Thm.2.1 \& Thm.2.2]{neel} using  purely  stochastic methods. He also studied the intersection problem in a class of minimal surfaces defined by stochastic properties \cite[Thm.5.1]{neel}.  In the  spirit of \cite{neel}, we  prove,  using potential theory tools, a strong half-space theorem between the classes of complete minimal surfaces with bounded curvature and  of  recurrent minimal surfaces  of $\mathbb{R}_{\raisepunct{.}}^{3}$ 


\begin{theorem}\label{main_thm_0}
Let $M$ be a recurrent minimal surface and $N$ be a complete  minimal surface with bounded curvature, immersed into $\mathbb{R}^{3}_{\raisepunct{.}}$ Then  $M\cap N\neq \emptyset$   unless they are parallel flat surfaces.\footnote{$M$ could be a plane minus a set of capacity zero parallel to a plane $N$.}
\end{theorem} 


A Riemannian manifold  is said to be recurrent (parabolic) if the standard Brownian motion  visits any open set at arbitrary large moments of time with probability 1 and it is transient otherwise. It is known that the recurrence of a manifold, not necessarily geodesically complete, can be described in terms of various analytic, geometric and potential theoretic properties (see \cite{grigoryan,PRS-PAMS,prs-memoirs}), for instance, it is equivalent to the following Liouville property: any bounded solution of the subequation $\triangle u \geq 0$  is constant.


%

The class  of  recurrent immersed minimal surfaces of $\mathbb{R}^{3}$ is  large. It contains all   complete minimal immersions of the complex plane $\mathbb{C}$ into $\mathbb{R}^{3}_{\raisepunct{,}}$ the complete  properly embedded minimal surfaces of $\mathbb{R}^{3}$ with finite genus \cite{meeks-perez-ros}, the complete minimal surfaces with quadratic volume growth, in particular, complete  surfaces with finite total curvature. In contrast,   the  first surface of Scherk is transient, see \cite{markvorsen,markvorsen-mcguinness-thomassen}.  The class of  recurrent  minimal surfaces is not contained in the class of complete properly immersed minimal surfaces nor on the class of complete minimal surfaces with bounded curvature. In the Apendix we  present examples of  recurrent non-proper minimal surfaces of  $\mathbb{R}^{3}$ with unbounded curvature. 
   
In \cite{RSS} Rosenberg, Schulze and Spruck, capturing the stochastic nature of  the Hoffmann-Meeks half-space theorem,  proved  a higher dimension half-space theorem for properly immersed minimal  hypersurfaces  of $M\times \mathbb{R}$, assuming that $M$ was a complete recurrent  $n$-manifold with bounded curvature. 

\begin{theorem}[Thm.1.2 of \cite{RSS}]Let $M$ be a complete recurrent  Riemannian $n$-manifold with bounded sectional curvature $\vert K_M\vert \leq \Lambda^2$ for some  $\Lambda \in \mathbb{R}$. Then any complete hypersurface minimally  and properly immersed in $M\times \mathbb{R}_{+}$ is a slice $M\times \{s\}$.\label{thmRSS}
\end{theorem} 

Recently, Theorem \ref{thmRSS} was extended, by Colombo, Magliaro, Mari and Rigoli \cite{cmmr}, to complete recurrent  Riemannian $n$-manifold with Ricci curvature bounded from below ${\rm Ric}\geq -(n-1) \Lambda^2_{\raisepunct{.}}$ Our second result is a version of \cite[Thm.1.2]{RSS} and  \cite[Thm.2]{cmmr} for minimal complete hypersurfaces with bounded curvature.


\begin{theorem}\label{main_thN}
Let $M$ be a complete recurrent Riemannian $n$-manifold with non negative Ricci curvature  $\text{Ric} \geq 0$, and  sectional curvature  bounded from above $K_M\leq \Lambda^2_{\raisepunct{.}}$ Then, any complete  hypersurface  $N$ minimally immersed in $M\times \mathbb{R}_{+}$ with bounded sectional curvature, equals a slice $M\times \{s\}$. 
\end{theorem}

\begin{remark}
Few  half-space theorems have been proved in others homogeneous $3$-spaces.  For instance, B. Daniel, W. Meeks, and H. Rosenberg \cite{daniel-meeks-rosenberg}, proved half-space theorems for properly immersed minimal surfaces of $Nil_{3}$ and $Sol_3$ where the half-space is defined by some distinguished minimal surfaces of these spaces, see also $\cite{daniel-hauswirth,menezes}$.
\end{remark}

The intersection problem for surfaces of $\mathbb{R}^{3}$ with constant mean curvature $H>0$,  called $H$-surface for short, was addressed by Ros and Rosenberg in \cite{ros-rosenberg}. Recall that   properly embedded $H$-surfaces $N$ separate $\mathbb{R}^3$ into two connected components and the mean convex side is the connected component of $\mathbb{R}^{3}\setminus N$ towards which the mean curvature vector field  points to.

\begin{theorem}[Ros-Rosenberg] \label{ros-rosenberg-thm}
A properly embedded  $H$-surface $M$ of $\mathbb{R}^{3}$ can not lie in the mean convex side of another properly embedded $H$-surface $N$.
\end{theorem}
Some half-space theorems for embedded $H$-surfaces in various homogeneous $3$-spaces were  proved in the appropriate settings. For instance, Rodriguez and Rosenberg  proved a half-space theorem for properly embedded $1$-surfaces of $\mathbb{H}^{3}$ in \cite{rodriguez-rosenberg},  Hauswirth, Rosenberg and Spruck  \cite{hauswirth-rosenberg-spruck}, Earp and Nelli in \cite{earp-nelli}, proved half-space theorems for properly embedded $1/2$-surfaces of $\mathbb{H}^{2}\times \mathbb{R}$, see  \cite{mazet-wanderley} for graphs  with $0<H<1/2$ in $\mathbb{H}^{2}\times \mathbb{R}$.  Using a fairly general approach, Mazet \cite{mazet} unified the proof of various half-spaces theorems for generic $3$-spaces, in particular for Lie groups with left invariant metric, see also \cite{mazet2013}.

In the second part of this article we are going to consider  half-space theorems for $H$-surfaces of $\mathbb{R}^{3}$ in the same vein of Theorems \ref{main_thm_0} \& \ref{ros-rosenberg-thm}. Let $N$ be a complete oriented properly immersed surface of $\R^3$ and denote by $\xi$ its unit normal vector field along $N$ for which $\overrightarrow{H}_{\!\!_N} = H_{\!_N}\xi$, $H_{\!_N}\geq 0$. Let $W$ be a connected component of $\R^3\setminus N$. Following ideas from \cite[Sec.4]{mazet}, we shall say that the mean curvature vector field of $N$ at $z_0 \in N \cap \partial W$ points into $W$ if, for any sequence $y_n \in W$ with $y_n \rightarrow z_0$ we have $y_n = \exp_z(t\xi(z))$ for some $0<t<\varepsilon$ and $z$ in a neighborhood $V \subset N$ of $z_0$. 

\begin{definition}\label{def_well_oriented} Let $M$ be a surface immersed into $\R^3_{\raisepunct{.}}$ A properly immersed surface $N$ into $\R^3_{\raisepunct{,}}$ disjoint from $M$, is said to be well-oriented with respect to $M$ if, $M$ lies in a connected component of $\R^3\setminus N$ for which the mean curvature vector field of $N$ points to. 
\end{definition}

Every complete oriented properly immersed minimal surface $N$ of $\R^3$ can be seen as well-oriented with respect to any surface immersed in $\R^3\setminus N$. Our next result gives a stochastic version of Theorem \ref{ros-rosenberg-thm}.

 \begin{theorem}\label{main_thm_1}
Let $M$ be an immersed surface of $\R^3_{\raisepunct{,}}$ and  $N$  an  oriented complete surface  properly immersed in $\R^3$ with bounded curvature. Then, unless $M$ and $N$ are parallel flat minimal surfaces, the surface $N$ can not be well-oriented with respect to $M$ provided that either
\begin{itemize}
\item[$a)$] $M$ is recurrent with mean curvature $\sup_{\!_M} \vert H_{\!_M}\!\vert  \leq \inf_{\!_N} H_{\!_N}$ or
\item[$b)$] $M$ is stochastically complete with $\sup_{\!_M}\vert  H_{\!_M}\!\vert < \inf_{\!_N} H_{\!_N}$.
\end{itemize} 
\end{theorem}
\begin{remark} 
This result should be compared with  \cite[Thm.1]{lmmg} proved assuming a bounded bending assumption.
\end{remark}

A Riemannian manifold $M$ is stochastically complete if the  diffusion process associated to the Laplacian $\triangle$ satisfies the conservation property  \begin{equation}\label{eq_heat_kernel}
\int_{M}p(t, x,y)d\mu(y)=1,
\end{equation}
for some/every $x\!\in\! M $ and all $t>0$. Here $p \in C^{\infty}\left((0, +\infty)\times  M \times M\right)$ is the heat kernel of $M$. The equation \eqref{eq_heat_kernel} has the following stochastic interpretation. The probability of the Brownian motion $X_t$ emanating  from $x$ to be found in $M$ is $1$, see  \cite{grigoryan, grigoryan_book}. The class of stochastically complete manifolds  contains  all complete manifolds with quadratic curvature decay or with quadratic exponential volume growth, as well as the properly immersed submanifolds of $\R^n$ with bounded mean curvature \cite{grigoryan,prs-memoirs}.

Among the  many equivalent characterizations for stochastic completeness,  we are going to use the following Liouville property: for all
$\lambda>0$, any bounded, non-negative solution of the subequation $\triangle u \geq \lambda u$  is identically zero. Therefore, recurrent manifolds are also stochastically complete. The converse statement is not true, for instance, the Scherk's first surface is stochastically complete and transient.

It is curious that the difference between the nature of intersection results for the class of minimal and $H$-surfaces is revealed by the threshold  $\lambda=0$ and $\lambda>0$ in the Liouville properties, translated as the conditions  a) and b). When $M$ is a stochastically complete surface of $\R^3$ satisfying $\sup_{\!_M}\vert H_{_M}\vert= \inf_{_N} H_{\!_N}>0$ we are able to prove only that ${\rm dist}(M,N)=0$. This result can be seen as version of \cite[Thm.5.1]{neel} for surfaces with positive mean curvature.


%
\begin{theorem}\label{main_thm_2} Let $M$ be a stochastically complete surface, and $N$ be a complete proper surface with bounded curvature, immersed in $\R^3_{\raisepunct{.}}$ If $N$ is well-oriented with respect to $M$ and $\sup_{\!_M}\vert H_{_M}\vert= \inf_{_N} H_{\!_N}>0$, then
$\dist(M,N) = 0$.
\end{theorem}


The Theorem \ref{main_thm_2} can be restated in terms of  relevant  geometric  conditions, sufficient for stochastic completeness as follow.

\begin{corollary}
Let $N$ be a complete embedded $H$-surface of $\R^3$ with bounded curvature and $M$ be an $H$-surface of $\R^3$  immersed in the mean convex side of $N$. Then,  $\dist(M,N)=0$ provided that either
\begin{enumerate}
\item[1.] $M$ is properly immersed,
\item[2.] $M$ has curvature $K_M(x)\geq -\rho^{2}(x)$, $\rho={\rm dist}_{M}(x_o,x)$, 
\item[3.] $M$ has volume growth ${\rm vol}(B_o(r)\cap M)\leq A\,e^{r^2}_{\raisepunct{,}}$ $A>0$.
\end{enumerate} 
\end{corollary}

Finally, we notice that the techniques developed to prove our results can be adapted to prove the Maximum Principle at Infinity between parabolic surfaces with non-empty boundary (possibly non-compact) and complete surfaces with bounded curvature immersed in $\R^3_{\raisepunct{.}}$ Roughly speaking, the Maximum Principle at Infinity generalizes Hopf's Maximum Principle for surfaces with constant mean curvature when the point of contact is at infinity. Several versions of the Maximum Principle at Infinity were proved in \cite{langevin-rosenberg,meeks-rosenberg-mp,soret-mp} for minimal surfaces, and in \cite{lima,meeks-lima} for $H$-surfaces, and later generalized in \cite{lmmg,meeks-rosenberg,ros-rosenberg}.


\begin{theorem}\label{main_thm_3}
Let $M$ and $N$ be two disjoint immersed surfaces of $\R^3_{\raisepunct{.}}$ Assume that $M$ is parabolic with non-empty boundary $\partial M,$ and $N$ is a complete surface with bounded curvature.
\begin{enumerate}
\item[1)] If both are minimal surfaces, then $\dist(M,N) = \dist(\partial M,N)$.
\item[2)]  If $\sup_{\!_M}\vert H_{\!_M}\vert \leq \inf_{_N} H_{\!_N}>0$, $N$ is proper and well-oriented with respect to $M$, then  $\dist(M,N) = \dist(\partial M,N)$.
\end{enumerate} 
\end{theorem}

From the stochastic viewpoint, a surface $M$ with boundary $\partial M$ is said to be parabolic if the absorbed Brownian motion is recurrent, that is, any Brownian path, starting from an interior point of $M$, reaches the boundary (and dies) in a finite time with probability $1$ (see \cite{perez-lopez}). From a potential-theoretic point of view (cf. \cite[Prop.10]{pessoa-pigola-setti}), it is equivalent to the following Ahlfors maximum principle: every bounded solution $u \in C^0(M)\cap W^{1,2}_{\text{loc}}(\inte M)$ of the subequation $\triangle u \geq 0$ on $\inte M$ must satisfies
\begin{eqnarray*}
\sup_{M} u = \sup_{\partial M} u.
\end{eqnarray*}
This notion of parabolicity for surfaces with boundary is weaker than the natural definition for which the Brownian motion reflects at $\partial M$ (see \cite{impera-pigola-setti,pessoa-pigola-setti}). 

\begin{remark}
A careful analysis on the proofs of our results suggests that they may hold for hypersurfaces on a general ambient space with curvature bounded from below, positive injectivity radius and non-negative Ricci curvature.
\end{remark}

%

\vspace{0.2cm}

\noindent \textbf{Acknowledgements.}\,We are grateful to  L. Mari for helpful  discussions  about regularity  issues regarding various parts of this manuscript. Special thanks to   Davi Maximo, the department of  Mathematics at University of Pennsylvania and Mathematisches Forschungsinstitut Oberwolfach, where part of this work was conducted, for their warm hospitality. This work was partially supported by CNPq-Brazil.
  

\section{Proof of Theorem \ref{main_thm_0}}\label{sec_minimal}
Suppose that $N$ has bounded curvature and assume  that $M\cap N=\emptyset$ and that $M$ and $N$ are not parallel flat surfaces.
We will split the proof  in two steps. In the first  we address the case where $N$ is embedded and in the second step we adapt the proof in the step 1 to  treat the general  case.

\vspace{2mm}

\noindent {\bf Step 1$\colon$}$N$ is embedded.
\vspace{2mm}

\noindent   If  $N$ is  embedded then there is a tubular $\varepsilon$-neighborhood $U(\varepsilon)=T_\varepsilon(N)$ which is  embedded for every $0<\varepsilon\leq 1/\vert \Lambda\vert$, where $K_{\!_{N}}\geq -\Lambda^2$ (see \cite[Thm.2]{soret}). For our purposes we are going to consider $0<\varepsilon<1/2\vert \Lambda\vert$. Since $N$ is two sided we can choose a smooth normal vector field $\eta$ to $N$ and decompose the tubular neighborhood $U(\varepsilon)=U_{-}(\varepsilon)\cup N \cup U_{+}(\varepsilon)$, where $\eta$ points toward the connected component $U_{+}(\varepsilon)$. Let   $t \colon  U(\varepsilon)\to \R$ be the signed distance function defined by
\begin{equation*}t(y)=\left\{
\begin{array}{lll}\,\,{\rm dist}(y, N) & {\rm if}& y\in U_{+}(\varepsilon),\\[0.2cm]
\!\!\!\!\!-{\rm dist}(y, N) & {\rm if}& y\in U_{-}(\varepsilon).
\end{array}\right.
\end{equation*}
Since ${\rm dist}_{\mathbb{R}^{3}}(M, N)=0$ (cf. \cite[Thm.5.1]{neel}),  we may assume that $M\cap U_{+}(\varepsilon)\neq \emptyset$.
\vspace{2mm}

Define $F\colon U_{+}(\varepsilon) \subset \mathbb{R}^3\to \mathbb{R}$ by
\begin{equation}\label{eqF}
F(y)= \log \left(\frac{2+\varepsilon\, c}{2+ 4\,c \,t(y)}\right),
\end{equation}
with $c \doteq \sup_{} \kappa >0$, where $\kappa \doteq \kappa_2 \geq 0$ is the non-negative principal curvature of $N$. Let $u \colon \varphi^{-1}(U_{+}(\varepsilon))  \to \mathbb{R}$ be given by $u=F\circ \varphi $, where $\varphi \colon M \to \mathbb{R}^{3}$ is the isometric minimal immersion of $M$ into $\mathbb{R}^{3}_{\raisepunct{.}}$ 
Clearly $u$ is smooth and  bounded. Since $M$ and $N$ are not parallel flat surfaces $M$ can not be a parallel surface $N^{t}$ to $N$ therefore $u$ is non-constant and $u\equiv 0$ on $\varphi^{-1}(\partial U_{+}(\varepsilon/4))$. We claim that 
\begin{eqnarray}\triangle_{_{M}}\! u \geq 0 & \text{on} &  \varphi^{-1}(U_{+}(\varepsilon/2)).\nonumber \end{eqnarray}
Indeed, consider the foliation $N^{t}$ by parallel surfaces to $N$ for $t \in (0, \varepsilon)$. For each $y\in N^{t}\cap U_{+}(\varepsilon/2)\cap M$ with coordinates $(x,t)\in N \times (0, \varepsilon/2)$  there is  an orthonormal basis $\{E_1,E_2\}\subset T_y N^{t}$ such that $\{E_1,E_2,\eta\}$ diagonalize the Hessian $\text{Hess}_{_{\R^3}}F$.
An easy computation yields
\begin{eqnarray*}
 \nabla_{_{\R^3}}F= - \frac{2c}{1+2c\, t}\, \eta & \text{and} &  \text{Hess}_{_{\R^3}}F = \frac{4c^2}{(1+2c\, t)^2} \nabla t \otimes \nabla t  - \frac{2c}{1+2c\, t} \nabla^2 t.
\end{eqnarray*}
Then, with respect to the splitting $\mathbb{R}\nabla t \oplus T_y \m_t$, the eigenvalues of $\text{Hess}_{_{\mathbb{R}^{3}}} F$ are 
\begin{eqnarray*}
\mu_1 =  - \frac{2c}{1+2c\, t}\frac{\kappa}{1+t\kappa}\raisepunct{,} \quad \mu_2 = \frac{2c}{1+2c\, t}\frac{\kappa}{1-t\kappa}\raisepunct{,} \quad \text{and} \quad \mu_3 =  \frac{4c^2}{(1+2c\, t)^2}\cdot
\end{eqnarray*} 
Since $0< 2t \leq \varepsilon < 1/2 \vert \Lambda \vert \leq 1/2c $, the monotonicity $\mu_1 \leq \mu_2 < \mu_3$ holds. Therefore,  applying  \cite[Lem.2.3]{kn:J-T} for the $2$-dimensional subspace $W\doteq  T_{y}{M}$ of $\mathbb{R}^3$ we have 
\begin{eqnarray*}
\triangle_{_{M}}\! u &=& \text{Tr }\text{Hess}_{_{\R^3}}F_{|_W}\\[0.2cm]
&\geq & \mu_1+\mu_2 \\[0.2cm]
& = & \frac{2c}{1+2c\, t}\frac{2t\ \kappa^2}{1-t^2\kappa^2}\\[0.2cm]
&\geq & 0.
\end{eqnarray*} 
Observe that  $M\cap U_{+}(\varepsilon/4)\neq \emptyset$ and 
 \begin{eqnarray}u\vert \varphi^{-1}( U_{+}(\varepsilon/4))>0,&  u\vert \varphi^{-1}(\partial U_{+}(\varepsilon/4))\equiv 0,& u\vert \varphi^{-1}( U_{+}(\varepsilon/2)\setminus U_+(\varepsilon/4))<0.\nonumber
 \end{eqnarray}The  function  $\overline{u}\colon M \to \mathbb{R}$ given by  $\overline{u}=\max\{u,0\}$ is continuous,  bounded and (non-negative) subharmonic  in the sense of distributions on a recurrent manifold thus it is constant by the Liouville property, see \cite[Thm.5.1]{grigoryan}. A contradiction.
 
\vspace{3mm}

\noindent{\bf Step 2$\colon$}$N$ is immersed.
\vspace{2mm}

\noindent As we have seen in the embedded case, we  need to construct a bounded weak solution $ u \in C^0(M)\cap W^{1,2}_{loc}(M)$ to the equation $\triangle u\geq 0$ on a recurrent manifold $M$. In view of Step 1 it would be natural to consider $u = F\circ \varphi \colon \varphi^{-1}\left(U(\varepsilon/2)\right) \to \R$ with $F = g\circ t_N$ given as in \eqref{eqF}, where $g \in C^{\infty}(\R_+)$, $\varphi : M \to \R^3$ is the isometric immersion and $t_N$ is the distance function to $N$. However, the distance function to $N$ is only Lipschitz continuous in general. The non-smoothness may occur when the set of self intersections  $\Gamma \subset \R^3$ of $N$ is non-empty, and the function $P\colon  U(\varepsilon) \to \overline{N}$ given by 
\[P(y)=\{z\in \overline{N} \colon {\rm dist}_{_{\R^3}}(y,z)={\rm dist}_{_{\R^3}}(y, N)\}\]
is a multivalued function. We will approach the subequation $\triangle u\geq 0$  considering solutions in the barrier sense. Recall that a function $u$ is said to satisfy $\triangle u \geq 0$ at a point $q$ in the barrier sense if, for any $\delta > 0$, there exists a smooth function $\phi_\delta$ around $q$ such that
\[
\begin{array}{cc}
\left\lbrace
\begin{array}{rl}
\phi_\delta =  u & \text{at} \ \ q, \\[0.2cm]
\phi_\delta  \leq u & \text{near } q,
\end{array}
\right.
& \quad \text{and} \qquad \triangle  \phi_\delta (q) > -\delta .
\end{array}
\]

Notice that since $N$ is minimal and has bounded curvature for each $x \in \overline{N}$  there exists a complete minimal surface $x\in L_x \subset  \overline{N}$ with  curvature $K_{L_x}\geq - \Lambda^2$. Thus, for a fixed $y\in U(\varepsilon/2)$ there exists, for each $z\in P(y)$, a simply connected, locally embedded neighborhood $V_z\subset L_{z}\subset \overline{N}$ of $z$  that is  graph over an open ball  $W_z\subset T_{z}L_{z}$ with radius uniformly bounded from below (see \cite{bessa-jorge-oliveira,rosenberg-hst}). Moreover, we can choose  $V_z$ so that ${\rm dist}_{\mathbb{R}^{3}}(y,z)\leq {\rm dist}_{\mathbb{R}^{3}}(y,z')$ for all $z'\in V_z$. Along each neighborhood $V_z$ we can consider a regular tubular neighborhood $C_z(\epsilon)=T_{\epsilon}(V_z)$ with radius $\epsilon$ and define the oriented distance function to $V_z$, $t_z  \colon C_z(\varepsilon) \to \R$, such that $t_z(y) > 0$. This yields  $C_z(\varepsilon)=C_z^{+}(\varepsilon)\cup V_z\cup C_z^{-}(\varepsilon)$. In order to construct a support function $\phi_\delta$ for the function $u$ we may select a neighborhood $V_z$, for $z \in P(y)$, and consider the function $F_z = g\circ t_{z}$. To avoid  the  analysis of the non focal points of the cut locus of the   boundary $\partial V_z$ of the surface $V_z$, we will introduce a supporting surface $S_z$ for $C_z^+(\varepsilon)$ at $z \in V_z$, that is, a smooth surface such that $z \in S_z$ and $C_z^+(\varepsilon)\cap S_z = \emptyset$. Following the agreement in \cite{lmmg}, modifying $S_z$ in a neighborhood around $z$, we may assume that $S_z$ is the boundary of a small, connected open set $B_{S_z} \subset C_z^-(\varepsilon)$. We can find a supporting surface $S_z^{\mu}$ for $C_z^+(\varepsilon)$ at $z \in V_z$ satisfying
$$H_{z}^{\mu}(z) > -\mu,$$
for any given $\mu > 0$, where $H_{z}^{\mu}$ is the mean curvature of $S_z^{\mu}$ (see \cite{lmmg}). This supporting surface can be constructed by deforming smoothly the boundary of a small ball $B \subset C_z^-(\varepsilon)$ touching $V_z$ at $z$. The following lemma says that we can choose $S_z^{\mu}$ so that $y \notin \text{cut}(S_{z}^{\mu})$.
 
\begin{lemma}[Lemma 1 in \cite{lmmg}]\label{lemma_lmmg} Fix $y \in U(\varepsilon/2)$ and a nearest point $z \in V_z$ to $y$. For a supporting surface $S_z$ at $z$, there exists $S'_{z}$, close to $S_z$ in the $C^\infty$ topology in a neighborhood of $z$, still supporting surface at $z$, and such that $y \notin \text{cut}(S'_{z})$.
\end{lemma}
Pick a point $ q \in M$ such that $y = \varphi(q) \in U(\varepsilon/2)$, $z \in P(y)$ and a neighborhood $V_z$ as described above. Given $\delta>0$, we will consider as a support function to $u = F\circ \varphi$ at $q$, the function $\phi_\delta \doteq F_{z}^{\mu}\circ \varphi$, where $F_{z}^{\mu} = g\circ t_{z}^{\mu}$, $t_{z}^{\mu}$ is the oriented distance function to $S_z^{\mu}$ with $t_{z}^{\mu}(y)>0$, and $S_z^{\mu}$ is a supporting surface as in Lemma \ref{lemma_lmmg}, for some $\mu = \mu(\delta)$ to be chosen later. Since $y \notin \text{cut}(S_z^{\mu})$ the support function $\phi_\delta$ is smooth in a small neighborhood of $q$. Furthermore, $\phi_\delta(q) = u(q)$, and taking a small ball $B_{\eta}(y)\subset U(\varepsilon/2)$ centred at $y$ and radius $\eta>0$, for every $\zeta\in B_{\eta}(y)$ it holds that 
\[
t(\zeta) \leq t_z^{\mu}(\zeta).
\]
Since $F$ is decreasing in $t$  we can assert that $\phi_\delta \leq u$ near to $q$.
In order to show that $u$ satisfies $\triangle_{_M} u (q) \geq 0$ in the barrier sense we are going to show that 
$$\triangle_{_M} \phi_\delta (q) > -\delta.$$  
Recall that 
\begin{equation*}
F_z^\mu(y)= \log \left(\frac{2+\varepsilon\, c}{2+ 4\,c \,t_{z}^{\mu}(y)}\right),
\end{equation*}
where now $c \doteq \max_{}\{\vert \kappa_1\vert,\vert \kappa_2\vert\} >0$, and $\kappa_1 \leq \kappa_2$ are the principal curvatures of $S_z^{\mu}$. Since $0< 2t < \varepsilon \leq 1/2c$, following up computations from Step $1$ we have
\begin{equation*}
 \text{Hess}_{_{\R^3}}F_z^\mu = \frac{4c^2}{(1+2c\, t_{z}^{\mu})^2} \nabla t_{z}^{\mu} \otimes \nabla t_{z}^{\mu}  - \frac{2c}{1+2c\, t_{z}^{\mu}} \nabla^2 t_{z}^{\mu},
\end{equation*}
whose eigenvalues are given by 
$$ \mu_1 =  \frac{2c}{1+2c\, t_{z}^{\mu}}\kappa_{1}^{t} , \quad \mu_2 = \frac{2c}{1+2c\, t_{z}^{\mu}}\kappa_{2}^{t}, \qquad \text{and} \quad \mu_3 =  \frac{4c^2}{(1+2c\, t_{z}^{\mu})^2}\raisepunct{,}$$
where
$$\kappa_1^t = \frac{\kappa_1}{1- t_{z}^{\mu}\kappa_1} \quad \text{and} \quad \kappa_2^t = \frac{\kappa_2}{1-t_{z}^{\mu}\kappa_2}$$ 
are the principal curvatures of the parallel surfaces to $S_z^{\mu}$ at $y$. We first observe that, independently of the sign of $\kappa_i$ $(i=1,2)$, it holds
\begin{eqnarray}\label{ineq_eigenv}
\mu_i \geq \frac{2c}{1+2c\, t_{z}^{\mu}}\kappa_{i}, \qquad \text{for $i = 1,2$}.
\end{eqnarray}
 The restriction on $\varepsilon\leq 1/2c$ and the inequality above  give us the monotonicity $\mu_1 \leq \mu_2 < \mu_3$. 
Again, applying  \cite[Lem.2.3]{kn:J-T} and inequality $H_{z}^{\mu} > -\mu$ we can write 
\begin{eqnarray}\label{ineq_laplacian_main}
\triangle_{_{M}} \phi_\delta &\geq & \mu_1+\mu_2 \nonumber \\[0.2cm] 
& = & \frac{2c}{1+2c\, t_{z}^{\mu}}H_{z}^{\mu}\\[0.2cm]
&\geq & -\frac{2c \mu}{1+2c\, t_{z}^{\mu}}.  \nonumber
\end{eqnarray} 
Then choose $\mu \doteq \delta/2c$ to conclude that $ \triangle \phi_\delta > -\delta.$

\vspace{2mm}

We can summarize the above discussion in the following proposition which could be of independent interest.

\begin{lemma}\label{mainlemma}
Let $N$ be a complete minimal surface immersed in $\R^3$ with bounded curvature. For any $\varphi \colon M \to \R^3$ complete minimal surface, there exist $\varepsilon>0 $ and a bounded function $u : M \to \R$ satisfying
\[ \triangle u  \geq 0 \quad \text{in the barrier sense,}\]
on the subset $\Omega_{\varepsilon} = \{ q \in M : 0< 2t(\varphi(q)) < \varepsilon\}$.
\end{lemma}

To conclude the proof in  the immersed case we apply the Lemma \ref{mainlemma} to show that $u=F\circ \varphi$ is a bounded,  subharmonic function in the barrier sense, therefore in the viscosity sense in $\varphi^{-1}(U(\varepsilon/2))$. Recalling that $u>0$ on $\varphi^{-1}(U(\varepsilon/4))$, and defining $\overline{u}=\max\{u, 0\}$ in $M$, by  \cite[Thm.1]{ishii} $\overline{u}$ is  a non-negative subharmonic in the sense of distributions and $\overline{u}\in C^0(M)\cap W_{loc}^{1,2}(M)$. We achieve the same contradiction as in the embedded case from the Liouville property \cite[Thm.5.1]{grigoryan}, see also \cite{mari_pessoa,mari_pessoa-2} for a direct proof from viscosity solutions.

\section{Proof of Theorem \ref{main_thN}}\label{sec_thm_1.3}

In the proof of Theorem \ref{main_thN} we intend to make explicit  how the geometry of the ambient space influences  this kind of intersection problem for minimal hypersurfaces. We will follow the same strategy applied in the proof of Theorem \ref{main_thm_0}. The main difference appears to be in the way to compute the Laplacian of the selected function, where we will use the ideas from \cite{mazet}. 

Recall that our assumptions on the sectional and Ricci curvatures of $M$ imply a uniform bound for the sectional curvature of the product ambient space $M\times \mathbb{R}$. Since $N$ is a minimal hypersurface with bounded sectional curvature the Gauss equation gives us a uniform bound for the second fundamental form of $N$. As a consequence of the extended Rauch's theorem (see \cite[Cor.4.2]{warner}) there exists a real value $\varepsilon > 0$ such that for every normal geodesic $\sigma$ issuing from a point $\sigma(0) \in N$ there is no focal points on $\sigma_{|_{[0,\varepsilon)}}$. By \cite[Prop.4.4]{docarmo} it means that the restriction exponential map $\text{exp}^\perp \colon (TN)^\perp \to M\times \mathbb{R}$ has no critical points in the tubular neighborhood $U(\varepsilon)$. Thus, for a fixed point in $U(\varepsilon)$ there is only one geodesic minimizing the distance to $N$.

Following ideas from Lemma \ref{mainlemma} we notice that for each point  $y \in U(\varepsilon)$, the projections points  that realize the distance from $y$ to $N$ are contained in $N$ or in the limit set of $N$. If a projection point $z$ lies on $N$ it is contained in a neighborhood $V_z \subset N$ which is locally embedded and has radius uniformly bounded from below. On the other hand, if the projection point $z$ lies on the limit set of $N$, then since our analysis is local, we can consider a uniform bound from below for the injective radius of the ambient space in order to guarantee the existence of a minimal neighborhood $V_z$, with radius uniformly bounded from below, contained in the limit set of $N$,  see the proof of  \cite[Thm.1.5]{bessa-jorge-oliveira} and \cite{bessa-jorge}. As in the proof of Theorem \ref{main_thm_0}, by Lemma \ref{lemma_lmmg}, at each projection point $z \in V_z$, for every $\mu > 0$, we can choose an embedding supporting surface $S_z^{\mu}$ at $z$ satisfying (see \cite[Lemm.1]{lmmg})
\begin{eqnarray}\label{modif_mu}
y \notin \text{cut}(S_{z}^{\mu}) \quad \text{and} \quad H_{z}^{\mu}(z) > -\mu.
\end{eqnarray}

Let us consider a slice $M\times \{s\}$, still called $M$, such that $\dist(M,N)=0$. Define a function $F \colon U(\varepsilon) \to \mathbb{R}$ given by $F = g\circ t_{_N}$ where 
\begin{equation*}
g(t)= \log \left(\frac{2+\varepsilon\, c}{2+ 4\,c \,t}\right),
\end{equation*}
$t_{_N} \colon U(\varepsilon) \to \mathbb{R}$ denotes the distance function to $N$, $c \doteq \max_{}\{\vert \kappa_1^t\vert,\ldots,\vert \kappa_n^t\vert\} >0$, and $\kappa_1^t \leq \cdots \leq \kappa_n^t$ are the principal curvatures of the parallel surfaces to $S_z^{\mu}$ for $0\leq t \leq\varepsilon$. We notice that $0<c<+\infty$ because the curvature bounds. Consider the function $u \colon \varphi^{-1}(U(\varepsilon/2)) \to \mathbb{R}$ given by $u = F\circ \imath$, with $\imath \colon M \to M\times \mathbb{R}$ be the inclusion isometric immersion. Again, our main claim is  that $u$ is a subharmonic function in the barrier sense on the subset $\varphi^{-1}(U(\varepsilon/2))$.


Given a point $ q \in M$ such that $y = \imath(q) \in U(\varepsilon/2)$, we assume the orientation $t_{_N}(y) >0$, and select a projection point $z \in \overline{N}$ for which a neighborhood $V_z$ has been associated. For any $\delta>0$ given, we will take $\mu = \mu(\delta)>0$ to be chosen later, and a supporting surface $S_z^{\mu}$ as in Lemma \ref{lemma_lmmg} satisfying \eqref{modif_mu}. Since $z \in S_z^{\mu}$, the oriented distance function to $S_z^{\mu}$, namely $t_{z}^{\mu}$, is smooth, touches $t_N$ from above on a neighborhood of $y$ and coincide with it at $y$. Thus, a support function to $u = F\circ \imath$ at $q$ is given by $\phi_\delta \doteq F_{z}^{\mu}\circ \imath$, where $F_{z}^{\mu} = g\circ t_{z}^{\mu}$.

It remains to compute the Laplace operator of $\phi_\delta$. 
For this, we take along $M$ an orthonormal basis $\{e_1,\ldots,e_n,e_{n+1}\}$ of $TM\times \mathbb{R}_+$ such that $e_{n+1} = \frac{\partial}{\partial t}$. Using this basis, and since $M$ is totally geodesic in $M\times \mathbb{R}_+$, we can write
\begin{equation}\label{eq_lapl_u_slice}
\triangle_{_M} \phi_\delta =  g'(t_{z}^{\mu})\sum_{i=1}^{n} Hess_{_{M\times \mathbb{R}_+}} t_{z}^{\mu}(e_i,e_i) + g''(t_{z}^{\mu}) \sum_{i=1}^{n} \langle \nabla_{_{M\times \mathbb{R}_+}} t_{z}^{\mu}, e_i\rangle^2.
\end{equation}

Let $S_t$ be the parallel surfaces to $S_z^\mu$, given by the image of a exponential map at time $t$, and by $\kappa_1^t, \ldots, \kappa_n^t$ its associated principal curvatures. Let $\{a_1,\ldots,a_n\}$ be an orthonormal basis of $TS_t$ which diagonalize the shape operator of $S_t$, and set $a_{n+1} = \eta$ where $\eta$ is the normal vector field along $S_t$ pointing towards $M$. The matrix of change of bases from $e_i$ to $a_i$ have the elements $(\lambda_{ij})_{1\leq i,j\leq n+1}$ defined by 
\[
e_i = \sum_{j=1}^{n+1} \lambda_{ij} a_j .
\]
With the above notation we can rewrite \eqref{eq_lapl_u_slice} as
\begin{equation*}
\triangle_{_M} \phi_\delta = g'(t_{z}^{\mu})\left(-\kappa_1^t(1- \lambda_{n+1,1}^2) - \cdots -\kappa_n^t(1 - \lambda_{n+1,n}^2 \right) + g''(t_{z}^{\mu})\left(1 - \lambda_{n+1,n+1}^2\right).
\end{equation*}

A main fact used in the proof of the previous theorem is the monotonicity of the mean curvature of the parallel surfaces along normal geodesics issuing from $S_z^\mu$, see \eqref{ineq_eigenv}. A sufficient condition for this monotonicity to hold is given by non-negativeness of the Ricci curvature of the ambient space, which in our case is guaranteed by the hypothesis ${\rm Ric}_M \geq 0$ (cf. \cite[Cor.3.5]{gray}). Therefore, using this monotonicity property and recalling that $g'(t_{z}^{\mu})<0$ we can estimate
\begin{eqnarray*}
\triangle_{_M} \phi_\delta &\geq & g'(t_{z}^{\mu})\left(\mu + \kappa_1^t \lambda_{n+1,1}^2 + \cdots + \kappa_n^t\lambda_{n+1,n}^2 \right) + g''(t_{z}^{\mu})\left(1 - \lambda_{n+1,n+1}^2\right).
\end{eqnarray*}
The $n$-dimensional spherical coordinates $(\theta_1,\theta_2,\ldots,\theta_n) \in [0,2\pi]\times [0,\pi]\times [0,\pi]$ can represent the unitary vector
$(\lambda_{n+1,1},\lambda_{n+1,2},\ldots,\lambda_{n+1,n+1})$ by
\begin{eqnarray*}
\lambda_{n+1,n+1} &=& \cos \theta_1 \\[0.1cm]
\lambda_{n+1,n} &=& \sin \theta_1 \cos \theta_2 \\[0.1cm]
\lambda_{n+1,n-1} &=& \sin \theta_1 \sin \theta_2 \cos \theta_3 \\[0.1cm]
& \vdots & \\[0.1cm]
\lambda_{n+1,2} &=& \sin \theta_1 \cdots \sin \theta_{n-1} \cos \theta_{n} \\[0.1cm]
\lambda_{n+1,1} &=& \sin \theta_1 \cdots \sin \theta_{n-1} \sin \theta_{n}.
\end{eqnarray*}
Applying these coordinates in the above estimate for the Laplacian together with the definition of the constant $0<c<+\infty$ we get
\begin{eqnarray*}
\triangle_{_M} \phi_\delta &\geq & g'(t_{z}^{\mu})\,\left(\mu + \kappa_1^t \left(\sin \theta_1 \cdots \sin \theta_{n-1} \sin \theta_{n}\right)^2 + \cdots + \kappa_n^t \left(\sin \theta_1 \cos \theta_2\right)^2 \right)\\[0.2cm]
& & + \ g''(t_{z}^{\mu})\left(1- \cos^2 \theta_1\right) \\[0.2cm]
&\geq & c\,g'(t_{z}^{\mu})\left(\left(\sin \theta_2 \cdots \sin \theta_{n-1} \sin \theta_{n}\right)^2 + \cdots + \left(\sin \theta_2 \cos \theta_2\right)^2 \right)\sin^2 \theta_1 \\[0.2cm]
& & + \ \mu g'(t_{z}^{\mu}) + g''(t_{z}^{\mu})\sin^2 \theta_1 \\[0.2cm]
&\geq & \mu g'(t_{z}^{\mu}) + \left(g''(t_{z}^{\mu}) + c\, g'(t_{z}^{\mu})\right)\sin^2 \theta_1 \quad \text{on } \ \varphi^{-1}(U(\varepsilon/2)).
\end{eqnarray*}
As before, reducing $\varepsilon$ if necessary, we can take $0<2t < \varepsilon \leq 1/2c$ and $\mu = \delta/2$. These choices lead us to conclude that 
$$\triangle_M \phi_\delta > -\delta.$$

Therefore, once established that $u$ is a bounded subharmonic function in the barrier sense on the subset $\varphi^{-1}(U(\varepsilon/2))$, we proceed  defining on $M$ the function $\bar{u} = \max\{u,0\}$ which will satisfy $\bar u \in C^0(M)\cap W^{1,2}_{loc}(M)$ and $\triangle_{M} \bar{u} \geq 0$ in the weak sense. Again, $\bar{u}$ will contradict Liouville property.


 

\section{Proof of Theorems \ref{main_thm_1} and \ref{main_thm_2}}\label{sec_cmc}


The proof of Theorems \ref{main_thm_1} and \ref{main_thm_2} follow the same strategy used in the proof  of Theorem \ref{main_thm_0} in the immersed case. Let us suppose, by contradiction, that $N$ is well-oriented with respect to $M$ and, up to an isomety of $\R^3$, we may assume that $\rm{dist}(M,N)=0$. Definition \ref{def_well_oriented} says that $M$ lies in an open connected component $W$ of $\R^3\setminus N$ for which the mean curvature vector field $\overrightarrow{H}_{\!_N} = H_{\!_N}\xi$ along $\partial W \subset N$ points into $W$. We recall that the boundary of $W$ is given as a union of smooth pieces of $N$ with non-negative mean curvature $H_{\!_N}$, and whose inner angles are not bigger than $\pi$ along an intersection set $\Gamma$. 

Similarly to the minimal case, there exists a regular tubular neighborhood $U_{+}(\varepsilon) \subset W$ with uniform radius $\varepsilon >0$ depending on the lower bound for the curvature of $N$. Let $t_{\!_N} : U_{+}(\varepsilon) \to \R$ be the distance function to $N$ which is a positive Lipschitz function. For any point $y \in U_{+}(\varepsilon)$ it is not hard to see that the nearest points to $y$ on $\partial W$ can not be on the part of $\Gamma$ where the inner angle is less than to $\pi$, otherwise the minimizing segment connecting $y$ to $\partial W$ will be normal to two different tangent planes. Therefore, for any point $z \in \partial W$ that minimizes the distance to $y \in U_{+}(\varepsilon)$, and any $\mu > 0$, we can deform one smooth piece of $N$ passing through $z$ to obtain a smooth supporting surface $S_z^{\mu}$ for $U_{+}(\varepsilon)$ at $z \in \partial W$ with mean curvature $H_{z}^{\mu}(z) > H_{\!_N}(z)-\mu$. Moreover, using Lemma \ref{lemma_lmmg} we can assume the oriented distance function to $S_z^{\mu}$, here called $t_{z}$, is smooth around $y$ and touches $t_{\!_N}$ from above at $y$.
 
Again, we set $c >0$ be the maximum norm of the principal curvatures of $N$ and consider the function $F \colon U_{+}(\varepsilon) \subset \mathbb{R}^3\to \mathbb{R}$ defined as $F = g\circ t_{\!_N}$, where
\begin{equation*}
g(t)= \log \left(\frac{2+\varepsilon\, c}{2+ 4\,c \,t}\right).
\end{equation*}
We also follow the convention $0<\varepsilon \leq 1/2c$. Set $\varphi \colon M \to \mathbb{R}^3$ be the isometric immersion of $M$, and define the function  $u = F\circ \varphi$ on $\varphi^{-1}(U_{+}(\varepsilon/2))$. We are going to prove that $u$ is a solution, in the barrier sense, of the subequation
\begin{eqnarray}\label{eq_u_barrier}
\triangle_{\!_M} u  \geq \frac{2c}{1+2c\, t_{z}}\left(\inf_{N} H_{\!_N} - \sup_{M} \vert H_{_{M}}\vert \right) \ \ \text{on} \ \ \varphi^{-1}(U_{+}(\varepsilon/2)) .
\end{eqnarray}
For any $x \in \varphi^{-1}(U_{+}(\varepsilon/2))$ and $\delta>0$, let us consider $\phi_\delta \doteq F_{z}^{\delta}\circ \varphi$, where $F_{z}^{\delta} = g\circ t_{z}$ and $t_{z}$ is the oriented distance function to $S_z^{\delta/2c}$ with $\varphi(x) \notin \text{cut}(S_z^{\delta/2c})$. Then, the function $\phi_\delta$ is a test function for $u$ at $x$. Arguing along similar lines from the proof of Theorem \ref{main_thm_0} we see that 
\begin{equation*}
 \text{Hess}_{_{\R^3}}F_z^\delta = \frac{4c^2}{(1+2c\, t_{z})^2} \nabla t_{z} \otimes \nabla t_{z}  - \frac{2c}{1+2c\, t_{z}} \nabla^2 t_{z},
\end{equation*}
whose eigenvalues are 
$$ \mu_1 =  \frac{2c}{1+2c\, t_{z}}\frac{\kappa_1}{1- t_{z}\kappa_1}\raisepunct{,} \quad \mu_2 = \frac{2c}{1+2c\, t_{z}}\frac{\kappa_2}{1-t_{z}\kappa_2}\raisepunct{,} \quad \text{and} \quad \mu_3 =  \frac{4c^2}{(1+2c\, t_{z})^2}\raisepunct{,}$$
where $\kappa_1 \leq \kappa_2$ are the ordered principal curvatures of $S_z^{\delta/2c}$. The monotonicity $\mu_1 \leq \mu_2 < \mu_3$ holds because $2\varepsilon\,c\leq 1$, as well as the inequality 
\begin{eqnarray*}
\mu_i =  \frac{2c}{1+2c\, t_{z}}\kappa_{i}^t \geq  \frac{2c}{1+2c\, t_{z}}\kappa_{i}, \quad \text{for $i = 1,2$}.
\end{eqnarray*}
Applying  \cite[Lem.2.3]{kn:J-T} we get
\begin{eqnarray*}
\triangle_{_{M}}\! \phi_\delta &=& \text{Tr}_{_{TM}}\!\text{Hess}_{_{\R^3}}F_z^\delta + \langle \nabla_{\!_{\R^3}}F_z^\delta, H_{\!_{M}} \rangle \\[0.2cm]
&\geq & \mu_1+\mu_2  - \frac{2c}{1+2c\,t_z}\sup_{M} \vert H_{\!_{M}}\vert\\[0.2cm] 
& \geq & \frac{2c}{1+2c\, t_{z}}\left(H_{\!_N} - \frac{\delta}{2c} - \sup_{M} \vert H_{\!_{M}}\vert \right)\\[0.2cm]
& \geq & \frac{2c}{1+2c\, t_{z}}\left(\inf_{N} H_{\!_N} - \sup_{M} \vert H_{\!_{M}}\vert \right) - \delta.  
\end{eqnarray*} 
Therefore, $u$ is a solution of the subequation \eqref{eq_u_barrier} in the barrier sense.

The proof of item $a)$ in Theorem \ref{main_thm_1} follows the same arguments employed on the previous proofs. For item $b)$ we notice that under the restriction on $\varepsilon$, setting $\lambda = \inf_{_N} H_{\!_N} - \sup_{_M} \vert H_{_{M}}\vert>0$ and using the inequality $s - 1 -\log s\geq 0$ for $s>0$, the function $u$ satisfies, in the barrier sense, 
\begin{eqnarray*}
\triangle_{_{M}}\! u &\geq & \lambda u \quad \text{on} \ \ \varphi^{-1}(U_+(\varepsilon/2)) .
\end{eqnarray*} 

Since $u$ vanishes only at $\varphi^{-1}(\partial U_+(\varepsilon/4))$ and it is subharmonic in the open set $\varphi^{-1}(U_+(\varepsilon/2))$, defining $\bar u \colon M \to \R$ as $\bar u = \max\{u,0\}$, we have that $\bar u$ is a bounded solution for the subequation $\triangle_{_M} \bar{u} \geq \lambda \bar{u}$ on $M$ in the weak sense, such that $\sup_M u > 0$. Again this contradicts the Liouville property for stochastic completeness \cite[Thm.5.1]{grigoryan}.

To prove Theorem \ref{main_thm_2} we shall assume by contradiction that $t_z \geq t_{\!_N} \geq 2\gamma$ for some constant $\gamma$ such that $0<2\gamma H< 1$, where $H = \sup_{_M} \vert H_{\!_M}\vert = \inf_{_N} H_{\!_N}$. We recall that the principal curvatures $\kappa_i^t$ of the parallel surfaces to $S_z^{\delta/2c}$ are given by
\begin{eqnarray*}
\kappa_i^t = \frac{\kappa_i}{1-\kappa_i t_z}\geq \kappa_i \quad \text{for} \ \ i=1,2.
\end{eqnarray*}
Denote by $H_z$ the mean curvature of $S_z^{\delta/2c}$ and assume $H_z > H - \delta/2c > H/2.$ Therefore, using that $\kappa_1 \leq H_{z}/2 \leq \kappa_2$, we have
\begin{eqnarray*}
\triangle_{_{M}}\! \phi_\delta 
&\geq & \frac{2c}{1+2c\,t_z}\left(\frac{\kappa_1}{1-\kappa_1 t_z} + \frac{\kappa_2}{1-\kappa_2 t_z} - H\right)\\[0.2cm] 
& \geq & \frac{2c}{1+2c\, t_z}\left(\kappa_1 + \frac{\kappa_2}{1-\gamma H_z}  - H_z -  \frac{\delta}{2c}\right)\\[0.2cm]
& \geq & \frac{2c}{1+2c\, t_z}\frac{\gamma H_z^2}{2-2\gamma H_z} - \frac{\delta}{1+2c\, t_z} \\[0.2cm]
& \geq & \frac{2c}{1+2c\, t_z}\frac{\gamma H^2}{4(2-\gamma H)} - \delta .
\end{eqnarray*} 
Taking $$\lambda = \frac{\gamma H^2}{4(2-\gamma H)}>0$$
we will conclude that $\triangle_{_M}\! u \geq \lambda u$ on $\varphi^{-1}(U_+(\varepsilon/2))$ in the barrier sense. The result can be finished by extending $u$ outside $\varphi^{-1}(U_+(\varepsilon/4))$ by zero and using the Liouville property for stochastic completeness \cite[Thm.5.1]{grigoryan}.

\section{Sketch of the proof of Theorem \ref{main_thm_3}}

As have seen before, up to  translation, we can assume that $\dist(M,N) = 0$. We just observe that the selected function $u$ used in the proof of all theorems is also a bounded solution of $\triangle u \geq 0$ on $\varphi^{-1}(U(\varepsilon/2)) \cap \inte M$, where $\varphi$ denotes the usual isometric immersion of $(M, \partial M)$ into $\R^3_{\raisepunct{.}}$ Furthermore, $\bar{u} = \max\{u,0\}$ belongs to $C^0(M)\cap W^{1,2}_{loc}(\inte M)$, and thus, it is a weak bounded  subharmonic function on $\inte M$. Since $M$ is assumed to be parabolic, the Ahlfors maximum principle \cite[Prop.10]{pessoa-pigola-setti} says that 
\begin{eqnarray*}
\sup_M u = \sup_{\partial M} u.
\end{eqnarray*}
To conclude, we only recall that $u(x) \to \sup_M u$ if and only if $\dist(\varphi(x),N) \to 0$.

%
%
%
%

\begin{appendices}

\section{}

In this appendix we provide two examples of recurrent minimal surfaces that are non-proper and have unbounded curvature. Following the ideas from \cite{andrade} we construct a non-trivial geodesically incomplete,  non-proper minimal immersion of $\mathbb{C}$ into $\mathbb{R}^{3}$ with unbounded curvature. The second example is a complete, non-proper minimal immersion of $\mathbb{C}$ into $\mathbb{R}^{3}$ with unbounded curvature.

\begin{example}
Consider an Enneper immersion $ \chi \colon \mathbb{C} \to \mathbb{C}\times \mathbb{R}$ given by
$$ \chi(z) = \left(L(z) - \overline{H}(z), h(z)\right),$$ where $L$ and $H$ are holomorphic functions defined by
\begin{eqnarray*}
L(z) = (r_1 - r_2)e^z,&  H(z) = -d\, e^{\left(\frac{r_1}{r_2}-1\right)z},
\end{eqnarray*}
and $h$ is a harmonic function defined as follow 
$$h(z) = - 4\left(\frac{d}{r_2}\right)^{\frac{1}{2}}\Big|\frac{r_2}{r_1}\Big|\vert r_1 - r_2\vert Re\left(i e^{\frac{r_1}{2r_2}z}\right).$$
We assume some non-degenerate assumptions for the parameters $r_1,r_2,d \in \mathbb{R}$, namely, $r_1 \neq r_2$ and $r_1 r_2 d \neq 0$, as well as some extra technical conditions
\begin{equation}
0 < r_1 < 4r_2 < 3r_1, \quad \frac{r_1}{r_2} \notin \mathbb{Q} \quad \text{and} \quad d = r_1 - r_2 > 0.
\end{equation}
The immersion $\chi(u+iv)$ is dense in an open subset of $\mathbb{R}^{3}_{\raisepunct{,}}$ its Gaussian curvature $K(u+iv) \to -\infty$ and $ds^2 = \lambda^2(u+iv)\vert dz\vert^2 \to 0$ as $u \to -\infty$. 
\end{example}

Our next example is a complete recurrent minimal surface which is neither properly immersed in $\R^3$ nor has bounded curvature.

\begin{example}
Let $f , g \colon \mathbb{C} \to \mathbb{C}$ be entire functions given by 
$$ f(z) = \frac{2}{\sqrt{\pi}} e^{r_{_1} z^2} \qquad \text{and} \qquad g(z) = e^{-r_{_2} z^2},$$
with the constants $r_{_1}, r_{_2}$ satisfying either $r_{_2} > r_{_1} > 0$ or $2r_{_2} > r_{_1} > r_{_2} > 0$. With the notations from \cite{barbosa_colares}, we 
define an immersion $\chi \colon \mathbb{C} \to \mathbb{R}^3$ by 
$$\chi(z) = (x_1(z),x_2(z),x_3(z))$$ where
\begin{eqnarray*}
x_1(z) &=& \text{Re}\int^z \frac{1}{2}\left(1-g^2\right)f\,dz = \frac{1}{\sqrt{\pi}}\text{Re}\int^z \left(e^{r_{_1} z^2} - e^{(r_{_1}-2r_{_2})z^2}\right)dz\\[0.2cm]
x_2(z) &=& \text{Re}\int^z \frac{i}{2}\left(1+g^2\right)f\,dz = -\frac{1}{\sqrt{\pi}}\text{Im}\int^z \left(e^{r_{_1} z^2} + e^{(r_{_1}-2r_{_2})z^2}\right)dz\\[0.2cm]
x_3(z) &=& \text{Re}\int^z g\,f\,dz = \frac{2}{\sqrt{\pi}}\text{Re}\int^z e^{(r_{_1}-r_{_2})z^2} dz.
\end{eqnarray*}
The above integrals defining coordinates $x_1$, $x_2$ and $x_3$ are given in terms of the error function and the imaginary error function. Explicitly, writing $\chi : \mathbb{C} \to \mathbb{C}\times \mathbb{R}$ we have
$$\chi(z) = \left(\frac{\text{erfi}(\sqrt{r_{_1}} \bar z)}{2\sqrt{r_{_1}}} -\frac{\text{erf}(\sqrt{2r_{_2} -r_{_1}} z)}{2\sqrt{2r_{_2} -r_{_1}}}\raisepunct{,} \frac{\text{Re}[\text{erf}(\sqrt{r_{_2} -r_{_1}} z)]}{\sqrt{r_{_2} -r_{_1}}}\right) \ \ \text{for } \ r_{_2} - r_{_1} > 0,$$
and
$$\chi(z) = \left(\frac{\text{erfi}(\sqrt{r_{_1}} \bar z)}{2\sqrt{r_{_1}}} -\frac{\text{erf}(\sqrt{2r_{_2} -r_{_1}} z)}{2\sqrt{2r_{_2} -r_{_1}}}\raisepunct{,} \frac{\text{Re}[\text{erfi}(\sqrt{r_{_1} -r_{_2}} z)]}{\sqrt{r_{_1} -r_{_2}}}\right) \ \ \text{for } \ r_{_1} - r_{_2} > 0.$$The  minimal surface $\chi$ is geodesically complete since its induced metric satisfies 
$$ds=\frac{1}{2}\vert f\vert\left(1+\vert g\vert\right)\vert dz \vert \geq  \frac{1}{\sqrt{\pi}}\vert dz\vert.$$  
The Gaussian curvature of $\chi$ is given by
\begin{eqnarray}
K(z) = - \left[\frac{4\sqrt{\pi}\,r_{_2} \vert z\vert}{\left(e^{\frac{r_{_1} +r_{_2}}{2}{\rm Re}(z^2)} + e^{-\frac{3r_{_2}-r_{_1}}{2}{\rm Re}(z^2)}\right)^2}\right]^2\cdot 
\end{eqnarray}
Let $z=te^{i\theta}$ be a point with polar coordinates $t>0$ and $\theta \in [0,2\pi)$. For $\varepsilon>0$ define  
$$C(\epsilon)= \bigcup_{i=1}^{4}\left(\frac{(2i-1)\pi}{4}-\epsilon, \frac{(2i-1)\pi}{4}+ \epsilon\right).$$
It is easy to see that if $\theta\not \in C(\epsilon)$, then $K(te^{i\theta})>-A\,t^2 e^{-B\,t^2},$ where $A= 16\pi r_{\!_2}^{2}$ and $B=(r_{\!_1}+r_{\!_2})\cos(\pi/2-2\epsilon)>0$. If $\theta \in C(\epsilon)$, then $K(te^{i\theta})\leq -A\,t^2 e^{-C\,t^2}$ where $C=(3r_{\!_2}-r_{\!_1})\cos(\pi/2-2\epsilon)>0$. In particular, 
$$ -At^{2}\geq \,K(te^{i\frac{k\pi}{4}}) \to -\infty , \quad \text{for} \ \ k=1,3,5,7.$$

Recall that the limit set $\lim \varphi $ of an immersion $ \varphi \colon M \to \mathbb{R}^{3}$ is defined as
$$\lim \varphi =\{q\in \mathbb{R}^{3}\colon \exists p_j\stackrel{in\,M}{\longrightarrow} \infty,\,\, {\rm  dist}_{\mathbb{R}^{3}}(\varphi(p_j),q)\to 0\,\, {\rm as}\,\, j\to \infty\}.$$
Along the curve $\gamma_{\theta}(t) = te^{i\theta}$ it is possible to  show that $\lim_{t\to \infty} \chi(\gamma_{\theta}(t))=q_{\theta} \in \R^3$  for  each $\theta \in  \{\frac{k\pi}{4}, k=1,3, 5,7\}.$ This proves that  the limit set $\lim \chi$ has at least the four points $\{q_{_{\pi/4}},q_{_{3\pi/4}},q_{_{5\pi/4}},q_{_{7\pi/4}}\}$. The points are these
\begin{itemize}
\item[] \hspace{-4mm}$q_{\pi/4}=(-\frac{1}{2\sqrt{2r_{\!_2}-r_{\!_1}}}\raisepunct{,} -\frac{1}{2\sqrt{r_{\!_1}}}\raisepunct{,} \frac{1}{\sqrt{r_{\!_2}-r_{\!_1}}}), \ $  $q_{3\pi/4}=(\frac{1}{2\sqrt{2r_{\!_2}-r_{\!_1}}}\raisepunct{,} -\frac{1}{2\sqrt{r_{\!_1}}}\raisepunct{,} -\frac{1}{\sqrt{r_{\!_2}-r_{\!_1}}}),$
\item[] \hspace{-5mm} $q_{5\pi/4}=(\frac{1}{2\sqrt{2r_{\!_2}-r_{\!_1}}}\raisepunct{,} \frac{1}{2\sqrt{r_{\!_1}}}\raisepunct{,}- \frac{1}{\sqrt{r_{\!_2}-r_{\!_1}}}), \ \,$
$q_{7\pi/4}=(-\frac{1}{2\sqrt{2r_{\!_2}-r_{\!_1}}}\raisepunct{,} \frac{1}{2\sqrt{r_{\!_1}}}\raisepunct{,} \frac{1}{\sqrt{r_{\!_2}-r_{\!_1}}})$.
\end{itemize}

Thus, $\chi\colon \mathbb{C}\to \mathbb{R}^{3}$ is neither proper nor has bounded curvature. To illustrate  see  the curve $x\to \chi (x+ix)$.
\begin{figure}[H]
\vspace{-0.3cm}
\begin{center}
\includegraphics[scale=0.3]{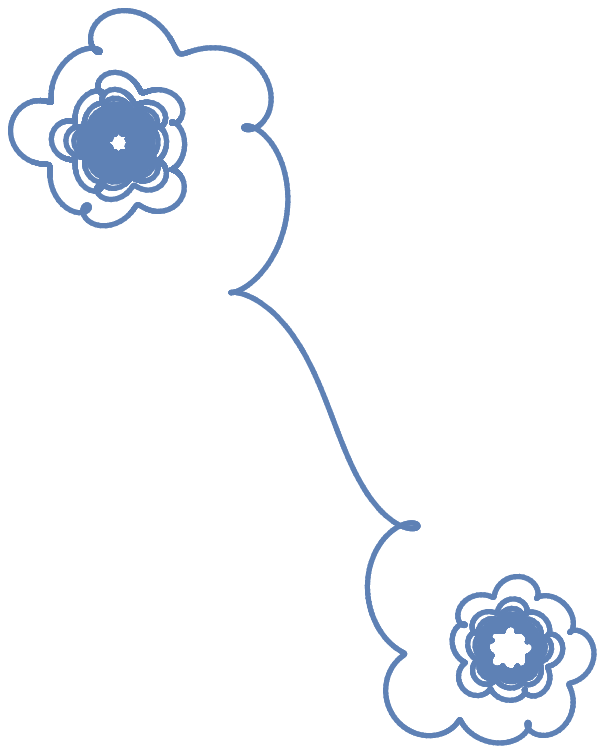}
\vspace{-0.6cm}
\caption{Curve $z = x+ix$ in a surface with $r_{_1} = 1, r_{_2} = 5$.}
\end{center}
\end{figure}

\end{example}
\end{appendices}

\vspace{0.3cm}
{\small
\begin{flushleft}
G. Pacelli Bessa \\
Departamento de Matem\'atica \\ Universidade Federal do Cear\'a\\
60455-760-Fortaleza, Brazil\\
bessa@mat.ufc.br.         
\end{flushleft}

\begin{flushleft}
Luquesio P. Jorge \\
Departamento de Matem\'atica\\ Universidade Federal do Cear\'a\\
60455-760-Fortaleza, Brazil\\
ljorge@mat.ufc.br. 
\end{flushleft}

\begin{flushleft}
Leandro F. Pessoa \\
Departamento de Matem\'{a}tica\\
Universidade Federal do Piau\'{i},\\
64049-550, Teresina - Piau\'i, Brazil\\
leandropessoa@ufpi.edu.br.
\end{flushleft}
}


\begin{thebibliography}{}

%
%
%
%
%
%
\bibitem{andrade}  P. Andrade, {\em A wild minimal plane in $\R^3$.} Proc. Amer. Math. Soc. \textbf{128}, n. 5 (2000) 1451--1457.
%
%
%
%
\bibitem{barbosa_colares} J. L. Barbosa, A. G. Colares, {\em Minimal Surfaces in $\mathbb{R}_{\raisepunct{.}}^3$} Lectures Notes in Mathematics \textbf{1195}, Springer-Verlag, Berlin, 1986.
%
%
\bibitem{bessa-jorge} G. P. Bessa, L. P. Jorge, {\em On properness of minimal surfaces with bounded curvature.} An. Acad. Brasil. Ci\^{e}nc. \textbf{75} (2003), n. 3, 279--284.

\bibitem{bessa-jorge-oliveira} G. P. Bessa, L. P. Jorge and G. Oliveira-Filho, {\em Half-space theorems for minimal surfaces with bounded curvature.} J. Differential Geom. \textbf{57} (2001), n. 3,  493--508.
%
%
%
%
%
%
%
\bibitem{cmmr} G. Colombo, M. Magliaro, L. Mari and M. Rigoli, {\em Bernstein and half-space properties for minimal graphs under Ricci lower bounds.} arXiv:1911.12054v2 [math.DG] 25 Jan 2020.

\bibitem{daniel-hauswirth} B. Daniel, L. Hauswirth, {\em Half-space theorem, embedded minimal annuli and minimal graphs in the Heisenberg group.} Proc. Lond. Math. Soc. (3) \textbf{98} (2009), n. 2, 445--470.
%
%
\bibitem{daniel-meeks-rosenberg} B. Daniel, W. H. Meeks III, H. Rosenberg, {\em Half-space theorems for minimal surfaces in ${\rm Nil}_3$ and ${\rm Sol}_3$}. J. Differential Geom. \textbf{88} (2011), n. 1, 41--59.

\bibitem{docarmo} M. P. do Carmo, {\em Riemannian Geometry}. Translated from the second Portuguese edition by Francis Flaherty. Mathematics: Theory $\&$ Applications. Birkhauser Boston, Inc., Boston, MA, 1992.


%
%
\bibitem{earp-nelli} R. Earp, B. Nelli, {\em A halfspace theorem for mean curvature $H = \frac{1}{2}$ surfaces in $\mathbb{H}^{2}\times \mathbb{R}$.} J. Math. Anal. Appl. \textbf{365} (2010), n. 1, 167--170.
%
%
%
%
%
%


\bibitem{lmmg} E. S. Gama, J. H. Lira, L. Mari, and A. A. Medeiros, {\em A barrier principle at infinity for varifolds with bounded mean curvature.} arXiv:2004.08946v1 [math.DG] 19 Apr. 2020.



\bibitem{gray} A. Gray, {\em Tubes.}  Second edition. With a preface by Vicente Miquel. Progress in Mathematics, 221. Birkh\"auser Verlag, Basel, 2004. 
%
%
\bibitem{grigoryan} A. Grigor'yan, {\em Analytic and geometric background of recurrence and non-explosion of the Brownian motion on Riemannian manifolds.} Bull. Amer. Math. Soc. (N.S.) \textbf{36} (1999), n. 2, 135--249.

\bibitem{grigoryan_book} A. Grigor'yan, {\em Heat kernel and analysis on manifolds.} AMS/IP Studies in Advanced Mathematics, \textbf{47}. American Mathematical Society, Providence, RI; International Press, Boston, MA, 2009.
%
%
%
%
%
\bibitem{hauswirth-rosenberg-spruck} L. Hauswirth, H. Rosenberg, J. Spruck, {\em On complete mean curvature $1/2$ surfaces in $\mathbb{H}^{2}\times \mathbb{R}$}. Comm. Anal. Geom. \textbf{16} (2008), n. 5, 989--1005.
%
\bibitem{hoffmann-meeks} D. Hoffmann, W. H. Meeks III , {\em The strong halfspace theorem for minimal surfaces.} Invent. Math. \textbf{101}, (1990), n. 2, 373--377.
%
%
%
\bibitem{impera-pigola-setti} D. Impera, S. Pigola, and A. G. Setti, {\em Potential theory for manifolds with boundary and applications to controlled mean curvature graphs.} J. Reine Angew. Math. \textbf{733} (2017), 121--159.

\bibitem{ishii} H. Ishii, {\em On the equivalence of two notions of weak solutions, viscosity solutions and distribution solutions.} Funkcial. Ekvac. \textbf{38} (1995), n. 1, 101--120.
%
%
%
%

\bibitem{kn:J-T} L. P. Jorge and F. Tomi, {\em The barrier principle for minimal submanifolds of arbitrary codimension.} Ann. Global Anal. Geom. \textbf{24} (2003), n. 3, 261--267.

%
%
\bibitem{jorge-xavier-ann} L. P. Jorge, F.  Xavier, {\em A complete minimal surface in $\mathbb{R}^{3}$ between two parallel planes.} Ann. of Math. $(2)$ \textbf{112} (1980), n. 1, 203--206.
%
%
%
%
%
%

\bibitem{langevin-rosenberg} R. Langevin and H. Rosenberg, {\em A maximum principle at infinity for minimal surfaces and applications.} Duke Math. J. \textbf{57} (1998), n. 3, 819--828. 


  
\bibitem{lima} R. F. Lima, {\em A maximum principle at infinity for surfaces with constant mean curvature in Euclidean space.} Ann. Global Anal. Geom. \textbf{20} (2001), n. 4, 325--343.


\bibitem{meeks-lima} R. F. Lima, W. H. Meeks III, {\em Maximum principles at infinity for surfaces of bounded mean curvature in $\mathbb{R}^3$ and $\mathbb{H}^3$.} Indiana Univ. Math. J. \textbf{53} (2004), n. 5, 1211--1223.



%
%
%
\bibitem{mari_pessoa} L. Mari, and L. F. Pessoa, {\em Duality between Ahlfors-Liouville and Khas'-minskii properties for non-linear equations.} 
Comm. Anal. Geom. \textbf{28} (2020), n. 2, 395--497.

\bibitem{mari_pessoa-2} L. Mari, and L. F. Pessoa, {\em Maximum principles at infinity and the Ahlfors-Khas'minskii duality: an overview.} Contemporary research in elliptic PDEs and related topics, 419--455, Springer INdAM Ser., \textbf{33}, Springer, Cham, 2019.
%
%
%
\bibitem{markvorsen}S. Markvorsen, {\em  Distance geometric analysis on manifolds. Global Riemannian geometry: curvature and topology}, 1--54, Adv. Courses Math. CRM Barcelona, Birkh\"{a}user, Basel, 2003. 
%

\bibitem{markvorsen-mcguinness-thomassen} S. Markvorsen, S. McGuinness, and C. Thomassen, {\em Transient random walks on graphs and metric spaces with applications to hyperbolic surfaces}, Proc. Lond. Math. Soc. \textbf{64} (1992), n. 1, 1--20.
%
%


\bibitem{mazet} L. Mazet, {\em A general halfspace theorem for constant mean curvature surfaces.} Amer. J. Math. \textbf{135} (2013), n. 3, 801--834.

%
\bibitem{mazet2013} L. Mazet, {\em The half space property for cmc $1/2$ graphs in $\mathbb{E}(-1, \tau)$.}  Calc. Var. Partial Differential Equations \textbf{52} (2015), n.3-4, 661--680.
%

\bibitem{mazet-wanderley} L. Mazet, G. A. Wanderley, {\em A half-space theorem for graphs of constant mean curvature $0<H<1/2$ in $\mathbb{H}^2\times \mathbb{R}$}. Illinois J. Math. textbf{59} (2015), n. 1, 43--53.


%
\bibitem{meeks-perez-ros} W. H. Meeks III, J. Perez and A. Ros, {\em The geometry of minimal surfaces of finite genus. II. Nonexistence of one limit end examples.} Invent. Math. \textbf{158} (2004), n. 2, 323--341.
%

\bibitem{meeks-rosenberg-mp} W. H. Meeks III, H. Rosenberg, {\em Maximum principles at infinity for minimal surfaces in flat three manifolds.} Comment. Math. Helv. \textbf{65} (1990), n. 2, 255--270. 


\bibitem{meeks-rosenberg} W. H. Meeks III, H. Rosenberg, {\em Maximum principles at infinity.} J. Differential Geom. \textbf{79} (2008), n. 1, 141--165. 
%



%
%
%
\bibitem{Meeks-Simon-Yau} W. H. Meeks III , L. Simon, S. T.  Yau, {\em  The existence of embedded minimal surfaces, exotic spheres and positive Ricci curvature.} Ann. Math. $(2)$ \textbf{116} (1982), n. 3, 621--659.
%




\bibitem{menezes} A. Menezes, {\em A half-space theorem for ideal Scherk graphs in $M\times \mathbb{R}$.} Michigan Math. J. \textbf{63} (2014), n. 4, 675--685.
\bibitem{neel}R. W. Neel, {\em A martingale approach to minimal surfaces.} J. Funct. Anal. \textbf{256} (2009), n. 8, 2440--2472.
%
%
%
%

\bibitem{perez-lopez} J. P\'erez, {\em Parabolicity and minimal surfaces.} Joint work with F. J. L\'opez. Clay Math. Proc., \textbf{2}, Global theory of minimal surfaces, 163--174, Amer. Math. Soc., Providence, RI, 2005.


\bibitem{pessoa-pigola-setti} L. F. Pessoa, S. Pigola, and A. G. Setti, \textit{Dirichlet parabolicity and $L^1$-Liouville property under localized geometric conditions.} J. Funct. Anal. \textbf{273} (2017), n. 2, 652--693.


\bibitem{PRS-PAMS} S. Pigola, M. Rigoli and A. G. Setti, \textit{A remark on the maximum principle and stochastic completeness.} Proc. Amer. Math. Soc. {\bf 131} (2003), n. 4, 1283--1288.
%
%
\bibitem{prs-memoirs}S. Pigola, M. Rigoli, A. G. Setti, \textit{Maximum principles on Riemannian manifolds and applications.} Mem. Amer. Math. Soc. \textbf{174} (2005), n. 822.
%
%
%
%
%
%
\bibitem{rodriguez-rosenberg} L. Rodriguez and H. Rosenberg, {\em Half-space theorems for mean curvature one surfaces in hyperbolic space.} Proc. Amer. Math. Soc. \textbf{126} (1998), n. 9, 2755--2762.

\bibitem{ros-rosenberg} A. Ros and H. Rosenberg, {\em Properly embedded surfaces with constant mean curvature.} Amer. J. Math. \textbf{132} (2010), n. 6, 1429--1443.

\bibitem{rosenberg-hst} H. Rosenberg,  {\em Intersection of minimal surfaces of bounded curvatures.} Bull. Sci. Math. \textbf{125} (2001), n. 2, 161--168.

%
%
\bibitem{RSS} H. Rosenberg, F. Schulze, J. Spruck, {\em The half-space property and entire positive minimal graphs in $M\times\mathbb{R}$}. J. Differential Geom. \textbf{95} (2013), n. 2, 321--336.

%
%
\bibitem{rosenberg-toubiana} H. Rosenberg, E. Toubiana {\em  A cylindrical type complete minimal surface in a slab of $\mathbb{R}^3$.} Bull. Sci. Math. \textbf{111} (1987), n. 3, 241--245.
%
%
%

\bibitem{soret-mp} M. Soret, {\em Maximum principle at infinity for complete minimal surfaces in flat $3$-manifolds.} Ann. Global Anal. Geom. \textbf{13} (1995), n. 2, 101--116.

\bibitem{soret} M. Soret, {\em Minimal surfaces with bounded curvature in Euclidean space.} Comm. Anal. Geom. \textbf{9} (2001), n. 5, 921--950.


%
%
%
%
%
\bibitem{warner} F. W. Warner, {\em  Extension of the Rauch comparison theorem to submanifolds.} Trans. Amer. Math. Soc. \textbf{122} (1966), 341--356.


%
\bibitem{xavier} F. Xavier, {\em  Convex Hull of Complete Minimal Surfaces.} Math. Ann. \textbf{269}  (1984), n. 2, 179--182.
%
%
%
%
%
%
%
%
%
%

\end{thebibliography}
\end{document}